\documentstyle[12pt]{article}

  \newcommand{\const}{\rm const}

  \newcommand{\argmax}{\rm  argmax}

\begin{document}

   \begin{center}

   {\bf Inversion of Tchebychev-Tchernov  inequality }\\

\vspace{4mm}

  {\bf Ostrovsky E., Sirota L.}\\

\vspace{4mm}

 Israel,  Bar-Ilan University, department of Mathematic and Statistics, 59200, \\

\vspace{4mm}

E-mails: eugostrovsky@list.ru \\
 \hspace{4mm}
sirota3@bezeqint.net \\

\vspace{5mm}

 {\bf Abstract} \\

\vspace{4mm}

 \end{center}

 \  We derive  in this article the  {\it lower}  bound for tail of distribution for the random variables (r.v.)  through a
 lower estimate for its moment generating functions (MGF).\par

\vspace{4mm}


\vspace{4mm}

\section{ Definitions.  Notations. Previous results.  Statement of problem.}

 \vspace{4mm}

  \ Let   $ \  (\Omega = \{ \ \omega \ \}, F, {\bf P})  \  $ be certain non-trivial Probability Space. Let also a numerical valued
 r.v. $ \ \eta = \eta(\omega) \ $ with correspondent  {\it tail function}

$$
T_{\eta}(x) \stackrel{def}{=} \max( {\bf P} (\eta \ge x),  \ {\bf P}(\eta < -x) \ ), \ x > 0, \eqno(1.0)
$$
 and a numerical valued non-negative function $ \ \nu = \nu(\lambda),  \ \lambda \in [1, b), \
b = \const \in (1, \infty], \ $ not necessary to be continuous or convex, be such that

$$
{\bf E} \exp(\lambda \eta ) \le \exp(\nu(\lambda)). \eqno(1.1)
$$
 \ Then

$$
T_{\eta}(x) \le \exp(- \nu^*(x)),  \ - \eqno(1.2)
$$
be the classical Tchebychev-Tchernov (Markov, Cramer) inequality. \par

 \ The function $ \ \lambda \to {\bf E} \exp(\lambda \ \eta) \ $ is named as a Moment Generating Function (MGF) for the r.v. $ \ \eta, \ $  write
$ \  \Theta_{\eta}(\lambda) = \Theta(\lambda) := {\bf E} \exp(\lambda \ \eta),  \ $
if of course it there exists in some non-trivial interval of the  values $ \ \lambda: |\lambda| < \lambda_0, \ \lambda_0 = \const \in (0, \infty], \ $
(Cramer's condition). It is closely related with the so - called Grand Lebesgue Spaces  [GLS], see e.g. [1], [3]-[6], [9]-[12] etc. It is also widely  used in
the theory of Great Deviations. \par

 \ Recall that the Cramer's condition for the r.v. $ \ \eta\ $  is quite equivalent to the following tail estimate

$$
\exists \ \mu = \const > 0 \ \Rightarrow  T_{\eta}(x) \le \exp(- \mu \ x), \ x \ge 0. \eqno(1.2a)
$$

\vspace{4mm}

 \ For instance, the function

$$
 \phi_{\eta}(\lambda):=  \ln {\bf E} \exp (\lambda \ \eta)
$$
 is named as the {\it natural } function for the r.v. $ \ \xi, \ $  is generating function for the correspondent GLS. \par
Another name: deviation function for the r.v. $ \ \eta, \ $ if of course this r.v. $ \ \eta \ $ satisfies the Cramer's condition. \par

 \ The complete description of such a functions $ \ \phi = \phi_{\eta}(\lambda) \ $ may be found in the book  [??OstrMono], pp. 22-24. \par

 \ Denote for simplicity

$$
T_{\eta}(x)  := \exp(- G(x)) = \exp(- G_{\eta}(x)), \eqno(1.3)
$$
so that

$$
 G_{\xi}(x) = |\ln T_{\eta}(x) \ |,
$$
 then  the estimate (1.2) may be rewritten as follows

$$
G(x) = G_{\eta}(x) \ge \nu^*(x). \eqno(1.4)
$$
 \ This function $ \  G(x) = G_{\eta}(x)  \  $ may be named as an {\it  exponential tail function } for the correspondent r.v. $ \ \eta. \ $  \par

\vspace{4mm}

  \ Here and in the sequel the transform $ \  \nu \to \nu^* \  $ will be denote an ordinary Young - Fenchel, or Legendre operator

$$
\nu^*(x) \stackrel{def}{=} \sup_{\lambda \in [1,b)} (\lambda \ x - \nu(\lambda)). \eqno(1.5)
$$

\vspace{4mm}

 \ {\bf To what extent can this inequality (1.2), or equally (1.4)  be reversed? }  \par

 \vspace{4mm}

 \ To be more precisely,   suppose   the  numerical valued  r.v. $ \ \xi = \xi(\omega) \ $ with correspondent  {\it tail function}
$ \ T_{\xi}(x)  = \exp(-G_{\xi}(x) \ $ and a numerical valued non - negative function $ \ \phi = \phi(\lambda),  \ \lambda \in [1, b), \
b = \const \in [1, \infty], \ $ not necessary to be continuous or convex, be such that

$$
{\bf E} \exp(\lambda \xi ) \ge \exp(\phi(\lambda)). \eqno(1.6)
$$

 \ Can we conclude (under some natural conditions) that

$$
\exists c = c(\phi) \in (0,1) \ \Rightarrow    T_{\xi}(x) \ge \exp(- \phi^*(c \ x)), \ x \ge x_0 = \const \in (0, \infty),  \eqno(1.7)
$$
or not? \par

\vspace{4mm}

 \ {\bf   Our claim in this article  is to  ground under some additional conditions the  inequality (1.7),
 i.e. to obtain the  opposite result for the inequality (1.2). } \par

\vspace{4mm}

\section{Auxiliary estimates from the theory of saddle-point  method.}

 \vspace{4mm}

 \ We must use in advance one interest and needed further integrals. Namely, let $  \  (X, M, \mu), \ X \subset R  \  $
 be non-trivial measurable space with non-trivial sigma finite measure  $ \ \mu. \ $ \par

 \ We assume at once  $ \ \mu(X) = \infty, \ $ as long as the opposite case is trivial for us.
We intend to estimate for sufficiently  greatest values of real parameter $ \ \lambda \ $,  say $  \ \lambda > e, \ $  the
following integral

$$
I(\lambda) := \int_X e^{  \lambda x -  \zeta(x)  } \ \mu(dx).  \eqno(2.1)
$$
assuming of course its convergence for all the sufficiently greatest values of the parameter $ \  \lambda. \ $ \par

 \ Here  $   \ \zeta = \zeta(x)   \  $ is non-negative measurable function, not necessary to be convex. \par

 \ We will use the main results obtained in the recent articles [10]-[12]. The offered below
estimates may be considered in turn as a some generalizations of the  saddle-point  method. \par

 \ We  represent now two methods for {\it upper} estimate $  I(\lambda) $ for sufficiently greatest values of the real parameter $ \ \lambda. \ $ \par

  \ Note first of all that if in contradiction the measure   $  \ \mu \ $ is finite: $ \  \mu(X) = M \in (0, \infty); \ $ then the
integral  $ \ I(\lambda) \ $ allows a  very simple estimate

$$
 I(\lambda)  \le M \cdot \sup_{x \in X}  \exp \left( \lambda x - \zeta(x)    \right) =
M \cdot \exp \left( \zeta^*(\lambda)  \right).  \eqno(2.2)
$$

 \ Let now $   \ \mu(X) = \infty $ and let $ \   \epsilon = \const \in (0,1); \ $
 let us introduce the following  auxiliary integral

$$
K(\epsilon) = K[\zeta](\epsilon) := \int_X e^{-  \epsilon \zeta(x) } \mu(dx). \eqno(2.3)
$$

 \ It will be presumed its finiteness at last for some positive value  $ \   \epsilon_0 \in (0,1); \ $
 then

 $$
   \forall \epsilon \ge \epsilon_0 \ \Rightarrow K(\epsilon) < \infty. \eqno(2.4)
 $$

 \ It is proved in particular in  [12] that

$$
I(\lambda) \le K(\epsilon) \cdot  \exp \left\{ (1 - \epsilon) \zeta^* \left(  \frac{\lambda}{1 - \epsilon} \right)  \right\}. \eqno(2.5)
$$

 \ As a slight consequence:

$$
I(\lambda) \le K(\epsilon) \cdot  \exp \left\{ \ \zeta^* \left(  \frac{\lambda}{1 - \epsilon} \right) \  \right\}, \eqno(2.6)
$$
and of course

$$
I(\lambda) \le \inf_{\epsilon \in (0,1)}
\left[ K(\epsilon) \cdot  \exp \left\{ (1 - \epsilon) \zeta^* \left(  \frac{\lambda}{1 - \epsilon} \right)  \right\} \right].
\eqno(2.7)
$$

 \ An opposite method, which was introduced in particular case in [10], [11], sections 1.2.
 Indeed,  let again  $ \   \epsilon = \const \in (0, 1).  \  $  Introduce a new function
$$
R(\epsilon) := \int_X e^{\zeta( ( 1 -\epsilon) \ x) - \zeta(x) } \ \mu(dx).  \eqno(2.8)
$$

 \ Then

$$
I(\lambda) \le R(\epsilon)  \ e^{  \zeta^*(\lambda/ (1 -\epsilon))  }.  \eqno(2.9)
$$

  \ Denote

$$
M[\zeta](\epsilon) := \min(K(\epsilon), \ R(\epsilon)), \ \epsilon \in (0,1). \eqno(2.10)
$$

 We obtained actually  the following compound estimate. \par

  \ {\bf Theorem 2.1.}  Suppose

$$
\exists \ \epsilon \in (0,1) \ \Rightarrow   M[\zeta](\epsilon) < \infty. \eqno(2.11)
$$
 \ Then

$$
I(\lambda) \le M[\zeta](\epsilon) \cdot   e^{  \zeta^*(\lambda/ (1 -\epsilon))  }, \ \lambda > 0. \eqno(2.12)
$$
 \ As a slight consequence:

$$
I(\lambda) \le  \inf_{\epsilon \in (0,1)} \left[ \ M[\zeta](\epsilon) \cdot   e^{  \zeta^*(\lambda/ (1 -\epsilon))  } \ \right], \ \lambda > 0. \eqno(2.12a)
$$

\vspace{4mm}

 \ {\bf Theorem 2.2.}  Suppose the non - negative r.v. $ \ \xi \ $ satisfies the Cramer's condition. Then

$$
 \forall \epsilon > 0 \ \Rightarrow K(\epsilon) := K(G[\xi], \ \epsilon) < \infty. \eqno(2.13)
$$

\vspace{4mm}

 \  {\bf Proof.} Denote for brevity  $ \   G(x) =  G[\xi](x), \ $  so that again

$$
T_{\xi}(x) = e^{-G(x)}, \ x \ge 0.
$$
 \ We have for the positive values $ \ \lambda > 0 \ $

$$
\infty >\lambda^{-1} \ \Theta_{\xi}(\lambda) = \lambda^{-1} \ {\bf E} \exp(\lambda \ \xi) = \int_0^{\infty} \exp( \lambda x - G(x) ) \ dx =
$$

$$
\sum_{n=0}^{\infty} \int_n^{n+1} \exp( \lambda x - G(x) ) \ dx \ge \sum_{n=0}^{\infty} \exp(\lambda n - G(n+1)),
$$
therefore

$$
G[\xi](n) \ge \lambda \ n - d(\lambda), \ \lambda > 0.
$$

 \ Further,

$$
K[G](\epsilon) = \int_0^{\infty} \exp(-\epsilon \ G(x)) \ dx = \sum_{n=0}^{\infty} \int_n^{n+1}   \exp(-\epsilon \ G(x)) \ dx \le
$$

$$
\sum_{n=0}^{\infty} \exp(-\epsilon \ G(n)) \le \sum_{n=0}^{\infty}  \exp ( - \epsilon (\lambda n - d(\lambda))) < \infty, \ \epsilon > 0,
$$
Q.E.D. \par

 \vspace{4mm}

\section{Main result: exponential lower bound for tail of distribution. Unilateral approach.}

 \vspace{4mm}

 \ So, let the inequality (1.6) be given. We need to introduce some another notations, definitions, and conditions. \par
 \ We can and will suppose without loss of generality that the r.v. $ \ \xi \ $ is non-negative, as long as we intend to
investigate the right-hand tail behavior. The case of opposite tail  may be studied analogously. \par
 \ Further, define the auxiliary function

$$
\phi_1(\lambda) = \phi_1[\phi](\lambda)  :=  \ln \left[ \ \frac{ \exp  \phi(\lambda) - 1}{\lambda} \ \right], \ \lambda \ge 1. \eqno(3.1)
$$

 \ Evidently  $  \  \phi(\lambda) > 1, \ \lambda \ge 1.  \  $ \par

\vspace{4mm}

\ {\bf Definition 3.1.} We will say that $ \ \phi \ - \ $ function belongs to the class $ \  W: \  \phi(\cdot) \in W, \ $ if the auxiliary function
$ \ \phi_1[\phi](\lambda)  \ $ is equivalent to the source one in the following sense

$$
\phi(\cdot) \in W \ \Leftrightarrow \exists c_1 = c_1[\phi] = \const \in (0,1], \ \phi_1(\lambda) \ge \phi(c_1 \ \lambda), \ \lambda \ge 1. \eqno(3.2)
$$
(The  inverse  inequality is trivial). \par

 \ This condition (3.2) is satisfied, if for example the function $ \ \phi(\cdot) \ $ is regular varying:

$$
\phi(\lambda) = \lambda^r \ L(\lambda),
$$
where $ \ r = \const > 0 \ $ and $ \  L(\cdot) \  $ is positive continuous slowly varying at infinity function. \par

\vspace{4mm}

\ {\bf Definition 3.2.}  The exponential tail function $ \ G = G(x) \ $ is said to be {\it super convex,}  write $ \  G = G(\cdot) \in SC,  \  $ iff it is convex
and  is also lower semi-continuous. \par

\vspace{4mm}

 \ {\bf Theorem 3.1.} \ {\bf A.} \ Suppose for some $ \  \epsilon \in (0,1) \ M[\phi_1](\epsilon) < \infty. \  $ Then

$$
\exp \left\{  -G^* \left( \frac{\lambda}{1 - \epsilon}    \right)   \right\} \le  \exp \left(-\phi_{2, \epsilon}(\lambda) \right), \eqno(3.3)
$$
where
$$
 \exp \left(-\phi_{2, \epsilon}(\lambda) \right) := \frac{1}{M[\phi_1](\epsilon)} \cdot \exp (-\phi_1(\lambda)). \eqno(3.3a)
$$

\vspace{4mm}

\ {\bf B.} If in addition $ \ \phi(\cdot)  \in W, \ $ then

$$
\exp \left\{  -G^* \left( \frac{\lambda}{1 - \epsilon}    \right)   \right\} \le  \exp \left(-\phi( c_2(\epsilon) \ \lambda) \right), \eqno(3.4)
$$

\vspace{4mm}

 \ {\bf C.} If in addition the exponential tail function $  \   G(x) \ $ is super convex: $ \ G (\cdot) \in SC, \ $ then

$$
G(x) \le \phi^*(c_3(\epsilon) \ x) ), \ x \ge 1, \ \eqno(3.5)
$$
or equally

$$
T_{\xi}(x) \ge \exp \left\{  \ -  \phi^*(c_3(\epsilon) \ x)  \ \right\}, \ x \ge 1. \eqno(3.5a)
$$

\vspace{4mm}

 \ {\bf Proof.} \par

\vspace{4mm}

 \ {\bf Preview.} We represent for beginning the sketch of proof, not to be  strict. Suppose in addition that the r.v. $ \ \xi \ $ has a logarithmical
convex density $ \ f(x); \ $ this means  by definition that

$$
f(x) = f_{\xi}(x) = e^{-v(x)},  x \ge x_0 = \const > 0,
$$
where  $ \  v = v(x) \  $ is convex lower semi-continuous function. We have

$$
e^{\phi(\lambda)} \le {\bf E} e^{\lambda \ \xi} = \int_R e^{\lambda x - v(x)} \ dx =: J(\lambda).
$$

 \ The last integral may be estimated as $ \ \lambda \to \infty \ $  by means of the saddle - point method

$$
\ln J(\lambda) \asymp \sup_x (\lambda \ x - \phi(\lambda)) = \phi^*(\lambda),
$$
and we conclude therefore that for all sufficiently great values $ \ \lambda, \ \lambda \ge \lambda_0 = \const \ge 1 \ $

$$
J(\lambda) \le e^{v^*(C_1 \ \lambda)},
$$
following

$$
v^*(C_2 \lambda) \ge \phi(\lambda).
$$
 \ Further,

$$
v^{**}(x) \le \phi^*(C_2 x), \ C_2 = C_2[\phi](\lambda_0) = \const \in [1, \infty).
$$
 \ By virtue of theorem Fenchel-Moreau $ \ v^{**}(x) = v(x), \ $ and we deduce   finally

 $$
 f(x) \ge \exp( - \phi^*(c_1 \ x) ), \ T_{\xi}(x) \ge \exp( - \phi^*(c_2 x) ),  \ x \ge x_0 = \const \ge 1.
 $$
 where $ \  c_1, c_2 = \const =  c_1, c_2[\phi](x_0) \in (0,1). \  $ \par

 \ Recall that always

$$
 T_{\xi}(x) \le \exp( - \phi^*(x) ),  \ x \ge x_0 = \const \ge 0.
$$

\vspace{4mm}

 \ Notice in addition that  $ \  \lim_{ x \to \infty} c_2(x_0) = 1. \ $  Thus, our estimations are essentially non-improvable as $ \ x \to \infty. $ \par
 \ The consequences of the last relation will be used further, in the seventh section. \par

\vspace{4mm}

 \  Actually, let's move on the complete proof.
 Denote by $ \  F(x) \ $ a function of distribution for the r.v. $ \ \xi: \ F(x) = {\bf P}(\xi < x), \ x \ge 0. \ $ We have
by means of integration ``by parts''

$$
{\bf E} e^{\lambda \ \xi }= \int_0^{\infty} e^{\lambda \ x} \ d F(x) = - \int_0^{\infty}  e^{\lambda \ x} \ d T_{\xi}(x) = {\bf P}(\xi \ge 0) +
$$

$$
 \lambda  \int_0^{\infty} e^{\lambda \ x} \ T(x) \ dx = 1 + \lambda  \int_0^{\infty} e^{\lambda \ x} \ T(x) \ dx
\ge e^{ \phi(\lambda)}, \ \lambda \ge 1,
$$
therefore

$$
\int_0^{\infty} e^{ \ \lambda x - G(x)  \ } \ dx  \ge e^{ \ \phi_1(\lambda)  \ }.  \eqno(3.6)
$$

 \ One can apply the inequality (2.12) of theorem 2.1:

$$
 M[G](\epsilon) \cdot   e^{  G^*(\lambda/ (1 -\epsilon))  }  \ge e^{\phi_1(\lambda)},  \ \lambda > 1, \eqno(3.7)
$$
which is equivalent to the propositions (3.3), (3.3a). The next assertion follows immediately from the direct definition of the set $ \ W: \ $

$$
G^* \left( \frac{\lambda}{1 - \epsilon}   \right) \ge  \phi( c_2(\epsilon) \ \lambda), \eqno(3.8)
$$
following

$$
G^{**} ( \lambda ) \le  \phi( c_4(\epsilon) \ \lambda). \eqno(3.9)
$$

 \  But  as long as  $ \ G \in SC, \ $ then $ \  G^{**} = G, \ $ (theorem of Fenchel-Moreau,) therefore

$$
G \left( \lambda \right) \le  \phi( c_4(\epsilon) \ \lambda),  \eqno(3.9)
$$
and finally
$$
T_{\xi}(x) \ge \exp \left\{  \ -  \phi^*(c_3(\epsilon) \ x)  \ \right\}, \ x \ge 1,
$$
Q.E.D. \par

 \vspace{4mm}

 \ {\bf Remark  3.1.} The estimate (3.5a) is trivially satisfied if the r.v. $ \ \xi \ $  does not satisfy the Cramer's condition. \par

\vspace{4mm}

\section{ Main result. Exponential level. Bilateral approach.}

 \vspace{4mm}

 \ Statement of the problem: given the {\it bilateral} inequality for some real valued random variable $ \ \xi: \ $

$$
\exp \left(  \phi_1(\lambda) \right) \le {\bf E} \exp(\lambda \ \xi) \le \exp \left(  \phi_2(\lambda) \right), \ |\lambda| < \lambda_0 = \const \in
(0, \infty], \eqno(4.1)
$$

$$
\phi_1, \ \phi_2 \in \Phi; \ \phi_1 (\lambda) \le \phi_2(\lambda).
$$

 \ Our purpose in this section is obtaining the lower estimate for the tail function for the r.v. $ \xi. \ $ \par

 \ We significantly weaken condition of the last section imposed on the function $ \ \phi_2(\cdot) \ $  and intend to obtain more qualitative
 up to multiplicative constant fine estimations. \par

\vspace{4mm}

 \ The investigated in this report problem  but for the random variables  satisfying the  Cramer's condition was
 considered, in particular, in the monograph [11], chapter 1, sections  1.3 , \ 1.4.  We intend to improve the obtained therein
results. \par

 \ We must introduce now some used notations. \\

\vspace{4mm}

$ \ \ \ \  \ \ \ \ \ \   S(\lambda, x) := \lambda \ x - \phi_2^*(x), \ \ \ \ \ \    S^*(\lambda) := \sup_x S(\lambda,x), \ $

\vspace{4mm}

$$
x_0 = x_0(\lambda) := \argmax_{x \in R} S(\lambda, \ x) =  \left[ \ \left(\phi_2^* \right)' \ \right]^{-1}(\lambda), \eqno(4.2)
$$

 \ Further, let $ \ x_- = x_-(\lambda) \ $  be arbitrary variable from the set $ \ (0, x_0(\lambda)), \ $ and $ \ x_+ = x_+(\lambda) \ $
be arbitrary variable from the set $ \ (x_0(\lambda)), \infty). \ $   The set of all this variables  $ \ x_-, x_+ \ $ will be denoted by
$ \  X_0 = X_0(\lambda). \ $\\

\ For instance, one can take

$$
x_-^{\Delta} =x_-^{\Delta}(\lambda)  := x_0(\lambda(1 - \Delta)), \  x_+^{\Delta} = x_+^{\Delta}(\lambda) := x_0(\lambda(1 + \Delta)),
\eqno(4.2a)
$$
$ \ \Delta = \const \in (0,1). \ $  More generally, one can choose also

$$
x_-^{\Delta(1)} =x_-^{\Delta(1)}(\lambda)  := x_0(\lambda(1 - \Delta(1))),
$$
$$
 x_+^{\Delta(2)} = x_+^{\Delta(2)}(\lambda) := x_0(\lambda(1 + \Delta(2))), \eqno(4.2b)
$$
$ \ \Delta(1), \Delta(2) = \const \in (0,1). \ $  \par

\vspace{4mm}

 \  Introduce also  the numerical valued function

$$
 z \to G_-(z) = G[x_-,x_+](z), \ z \ge e
$$
as follows:\\

$$
    G_-(x_-(\lambda)) = G_-[x_-(\cdot), x_+(\cdot)](x_-(\lambda))  = G_-[x_-, x_+](x_-(\lambda)) \stackrel{def}{=}
$$

$$
 e^{ - \lambda \ x_+(\lambda) } \times
 \left\{ \  e^{\phi_1(\lambda)} - \lambda \frac{e^{S(\lambda, x_-(\lambda))}}{S'_x(\lambda, x_-(\lambda))} -
 \lambda \frac{e^{S(\lambda, x_+(\lambda))}}{|S'_x(\lambda, x_+(\lambda))|} \right\}, \eqno(4.3)
$$
as well as its ``closure''

$$
\overline{G}_-(z) \stackrel{def}{=} \sup_{ (x_-,x_+) \in X(\lambda)  } G_-[x_-,x_+](z). \eqno(4.3a)
$$

\vspace{4mm}

\ {\bf Theorem 4.1.} We conclude under formulated above notations and restrictions

$$
T_{\xi}(z) \ge \overline{G}_-(z), \ z \ge 1. \eqno(4.4)
$$

\vspace{4mm}

{\bf Proof. } We have after integration by part assuming $ \ \lambda > 0 \ $

$$
\lambda^{-1} \ e^{\phi_1}(\lambda) \le \ \int_{-\infty}^{\infty} e^{\lambda x} \ T_{\xi}(x) \ dx = I_1 + I_2 + I_3, \eqno(4.5)
$$
where

$$
I_1 = I_1(\lambda) = \int_{-\infty}^{x_-} e^{\lambda x} \ T_{\xi}(x) \ dx \le \int_{-\infty}^{x_-} e^{\lambda x - \phi_2^*(x)} \ dx =
$$

$$
\int_{-\infty}^{x_-} e^{S(\lambda,x)} \ dx \le \int_{-\infty}^{x_-}\exp \left\{ S(\lambda,x_-) + S'_x(\lambda, x_- )(x - x_-) \ \right\} \ dx =
$$

$$
\frac{\exp  S(\lambda,x_-)}{S'_x(\lambda, x_- ) }.
$$
  \ We find analogously

$$
I_3 = I_3(\lambda) = \int^{\infty}_{x_+} e^{\lambda x} \ T_{\xi}(x) \ dx \le \int^{\infty}_{x_+} e^{\lambda x - \phi_2^*(x)} \ dx =
$$

$$
\int^{\infty}_{x_+} e^{S(\lambda,x)} \ dx \le \int^{-\infty}_{x_+}\exp \left\{ S(\lambda,x_+) -| S'_x(\lambda, x_+ )| \ (x - x_+) \ \right\} \ dx =
$$

$$
\frac{\exp  S(\lambda,x_+)}{|S'_x(\lambda, x_+ )| }.
$$

 \ Finally,

$$
I_2 = I_2(\lambda) \le  \lambda \int_{x_-}^{x_+} T(x_-)  \  e^{\lambda x} \ dx <  T(x_-) \ e^{\lambda x_+}. \eqno(4.6)
$$

 \  We deduce after substituting into (4.5) and solving obtained inequality relative the variable  $  \ T(x_-) \ $

$$
T_{\xi}(z) \ge {G}_-[x_-, x_+]_-(z), \ z \ge 1.
$$

 \ It remains to make the optimization subject to the natural limitation $ \  (x_-,x_+) \in X(\lambda): \  $

$$
T_{\xi}(z) \ge  \sup_{ (x_-,x_+) \in X(\lambda)} {G}_-[x_-,x_+](z) =\overline{G}_-(z), \ z \ge 1,
$$
\ Q.E.D. \\

 \ {\bf Remark 4.1.} Of course, one can take accept as  the values $ \ (x_-, x_+) \ $ the variable from (4.2b)

$$
x_- := x_0(\lambda(1 - \Delta(1))),  \  x_+ := x_0(\lambda(1 + \Delta(2))), 0 < \Delta(1), \Delta(2) < 1.\eqno(4.7)
$$

\vspace{4mm}

\section{Consequences. Particular cases.}

\vspace{4mm}

 \  It seems quite reasonable to choose  the values $ \  x_{\pm} =   x_{\pm}(\lambda)   \  $  for sufficiently greatest values $ \ \lambda,\ $
say by definiteness $ \ \lambda \ge e, \ $ as before

$$
 x_{\pm}(\lambda) := x_0(\lambda(1 \pm \Delta)), \ \Delta = \const \in [0, 1/2],
$$
or for simplicity

$$
 x_{\pm}(\lambda) := x_0(\lambda(1 + \Delta)), \ \Delta = \const \in [-1/2, 1/2]. \eqno(5.0)
$$

 \ We impose also the following natural condition on the function $ \ \phi = \phi(\lambda) := \phi_2(\lambda) \ $

$$
\inf_{ 0 \ne \Delta \in [-1/2, 1/2]} \ \inf_{\lambda \ge e}
\left[ \frac{ S(\lambda, x_0) -  S(\lambda, x_0(\lambda(1 + \Delta)))}{S(\lambda, x_0) \ \Delta^2} \right] =: V[\phi] = V > 0. \eqno(5.1)
$$

 \  We consider in this section only the case when the function $ \ \phi_1(\lambda) \ $ satisfies the restrictions of the form

$$
\phi(\lambda) \ge
\phi_1(\lambda) \ge (1 - \delta^2) \ \phi(\lambda), \ \lambda \ge e, \ \delta = \const \in (0, 1/2), \eqno(5.2)
$$
so that

$$
\exp \left( (1 - \delta^2) \  \phi(\lambda) \right) \le {\bf E} \exp(\lambda \ \xi) \le \exp \left(  \phi(\lambda) \right), \ |\lambda| < \lambda_0 = \const \in
(0, \infty]. \eqno(5.2a)
$$

 \ Further,  note that

$$
\lambda \ x_0(\lambda) - \phi(\lambda) = \phi^*(x_0(\lambda)), \ \lambda \ge 1.
$$

 \ Therefore, it is naturally to suppose in addition

$$
\exists \ c_0 \in (0,\infty) \ \forall \lambda \ge e, \ \delta \in (0,1/2) \ \Rightarrow
$$

$$
\lambda \ x_0(\lambda (1 + \delta)) - (1 - \delta^2) \phi(\lambda) \le (1+c_0 \delta) \phi^*(x_0(\lambda(1 - \delta))). \eqno(5.3)
$$

\  Both the conditions (3.1), (3.3)  are satisfied if for example the function $ \ \lambda \to \phi(\lambda) \ $ is  sufficiently smooth  regular varying:

$$
\phi(\lambda) =\phi[p; L](\lambda) :=
 p^{-1} \ |\lambda|^p \ L(\ |\lambda| \ ), \ |\lambda| \ge e, \eqno(5.4)
$$
where $ \ p = \const > 1, \ $ and $ \  L(\cdot) \  $ is some positive continuous differentiable slowly varying at infinity function.
 For instance,  one can take

$$
\phi(\lambda) = \phi_{2}(\lambda) := 0.5 \ \lambda^2, \ \lambda \in R, \ -
$$
the so - called subgaussian case; as well as many popular examples

$$
\phi(\lambda) =\phi[p]_r(\lambda) :=
 p^{-1} \ |\lambda|^p \  [\ln (e + |\lambda|)]^r, \ r = \const, \ |\lambda| \ge 1. \eqno(5.4a)
$$

\vspace{4mm}

 \ {\bf Theorem 5.1.} Suppose that all the formulated before conditions, (2.1),  (3.1), (3.2), (3.3) imposed on the r.v.  $ \ \xi \ $
are satisfied. Then for some finite positive constant $ \  c = c(\phi)  \in (0, 1/(2\delta)) \  $

$$
T_{\xi}(z) \ge \exp \left\{ \ - (1 - c \delta) \phi^*(z/(1 - c \delta)) \  \right\}, \ z \ge e.   \eqno(5.5)
$$

\vspace{4mm}

 \ {\bf Proof.} One need only to apply the assertion of theorem 2.1., choosing $ \ \Delta = C_2 \ \delta, \  \Delta  \in (0, 1/2). \ $
We omit some simple calculations. \par

\vspace{4mm}

 \ {\bf Remark 5.1.} Suppose that for certain function $ \  \phi(\cdot) \in \Phi  \ $ and for some random variable $ \ \xi \  $

$$
\exp \left( (1 - \delta^2) \phi(\lambda) \right) \le {\bf E} \exp(\lambda \ \xi) \le \exp \left(  \phi(\lambda) \right),
$$

$$
\forall \lambda: \ |\lambda| < \lambda_0 = \const \in (0, \infty], \ \delta = \const \in (0,1/2].  \eqno(5.6)
$$

 \ Since

$$
\left[ \  (1 - \delta^2) \ \phi \   \right]^*(x) = \sup_{\lambda} \left(  \ \lambda x - (1 - \delta^2) \  \phi(\lambda)  \right)=
$$

$$
(1 - \delta^2) \sup_{\lambda} \left( \frac{\lambda}{1-\delta^2} \ x - \phi(\lambda)   \right) =
(1 - \delta^2) \ \phi^* \left( \frac{x}{1 - \delta^2}  \right),
$$
 it is natural to wait (our hypothesis!) that

$$
T_{\xi}(z) \ge \exp \left\{ \ - (1 - c_2 \delta^2) \phi^*(z/(1 - c_2 \delta^2)) \  \right\}, \ z \ge e;
$$
but we have  grounded  (under formulated above conditions) only the more weak estimate (3.5). \par

\vspace{4mm}

\ {\bf Example 3.1.} Let $ \ \phi(\lambda)  = 0.5 \ \lambda^2, \ \lambda \in R,\ \delta = \const \in (0,1/2);  \ $  then all the assumptions formulated above
hold true. In detail: $ \ \delta = \const \in (0, 1/2), \ $ and the (meal zero) r.v. $ \ \xi \ $ is such that

$$
\exp \left( (1 - \delta^2) \  \lambda^2/2) \right) \le {\bf E} \exp(\lambda \ \xi) \le \exp \left(  \lambda^2/2 \right), \ \lambda \in R,  \eqno(5.7)
$$
 a  subgaussian case. We derive by virtue of theorems 2.1-3.1

$$
T_{\xi}(z) \ge \exp \left\{ - 0.5 \ z^2(1 + \ c \delta) \ \right\},  \ z \ge 0,
$$
but in accordance with our hypotheses

$$
T_{\xi}(z) \ge \exp \left\{ - 0.5 \ z^2(1 + \ c \delta^2) \ \right\}.
$$
 \  {\it This is an open problem.} \par

\vspace{4mm}

\section{Main result. Power level.}

 \vspace{4mm}

 \ We intend in this section to deduce the lower bound for the  tail of probability  for the r.v. $ \ \xi: \ \ T_{\xi} (x), \ x \ge 1 \ $  through
its moment estimates, or equally through its Lebesgue-Riesz $ \ L_p = L_p(\Omega) \ $ norm

$$
 |\xi|_p = \left[ {\bf E}|\xi|^p \right]^{1/p}, \ p \in [1,b), \ b = \const \in (1, \infty]. \eqno(6.0)
$$

 \ This case may be easily reduced to the considered one. More detail, one can suppose

\vspace{4mm}

$$
|\xi|_p \ge \zeta(p), \ 1 \le p < b, \ b = \const \in (1, \ \infty], \eqno(6.1a)
$$
 an unilateral inequality or

$$
\zeta(p) \le  |\xi|_p \le \psi(p), \ 1 \le p < b, \ b = \const \in (1, \ \infty], \eqno(6.1b)
$$
a bilateral estimate. \par

 \ The relation, e.g. (6.1b) may be rewritten as follows. Put $ \ |\xi| = \exp(\theta), \ $  so that

$$
T_{\xi}(x) = T_{\theta}(\ln x), \ x \ge e,
$$

 and denote

 $$
 \phi_1(\lambda) = \lambda \ \ln \zeta(\lambda), \ \phi_2(\lambda) = \lambda \ \ln \psi(\lambda), \ \lambda \in [1,b);
 $$
then

$$
e^{\phi(\lambda)} \le {\bf E} \  e^{\lambda \ \theta} \le e^{\phi_2(\lambda)}. \eqno(6.2)
$$
 \ It remains to apply the results of foregoing sections 2-5. \par

\vspace{4mm}

 \ {\bf Example 6.1.} Let the r.v. $ \ \eta \ $ be such that

$$
 \forall p \in [1,b) \ \Rightarrow |\eta|_p  \ge C \ (b-p)^{-\beta},
$$

$$
 b = \const \in (1, \infty), \ \beta = \const \in (0, \infty). \eqno(6.3)
$$

 \ The examples of these variables may be found, e.g. in [1], [11] and so one. We conclude from (6.3) after some calculations

$$
{\bf P}(|\xi| > x) \ge C(b, \beta) \ x^{- \gamma}, \ x \ge 1, \  \exists \gamma = \const \in (1, b). \eqno(6.4)
$$

\vspace{4mm}

 \ {\bf Example 6.2.} Suppose the r.v. $ \ X \ $ is such that

$$
|X|_p \asymp p^{1/m}, \ p \ge 1, \ m = \const \in (0, \infty), \eqno(6.5)
$$
the both extremal "boundary" cases $ \ m = 0 \ $ or $ \ m = \infty \ $  are trivial. \par

 \ Many practical examples of such a r.v. may be found in the articles [2], [7]-[8],
devoting to a mesoscopic physics. The authors obtained in particular the mild and great deviations for the sequence of these variables. \par

 \ For instance, there exists the r.v. $ \ Y \ $ for which

$$
{\bf E}|Y|^p = C \ D^p\ \frac{\prod_{j=1}^J \Gamma(a_j \ p + b_j)}{ \prod_{k=1}^K \Gamma(a_k' \ p + b_k')},
$$
 where $ \ \Gamma(\cdot) \ $ is ordinary Gamma-function, see [7]-[8]. \par

 \ We conclude that under relation (6.5)

$$
\exp(-C_2 x^m ) \le T_{X}(x) \le \exp(-C_1 x^m ), \ x \ge 1,
$$

$$
 C_1, C_2 = \const \in (0,\infty), C_1 \le C_2. \eqno(6.6)
$$

 \ Notice that the r.v. $ \ X \ $ from (6.5) satisfies the Cramer's condition only when $ \ m  \ge 1. \ $ \par

\vspace{4mm}

\section{ Tauberian theorem. }

 \vspace{4mm}

 \ {\it 1.  Ordinary approach.} \par

 \vspace{4mm}

 \ {\bf Theorem 7.1.} Assume the function $ \ \phi(\cdot) \ $ satisfies all the conditions of theorem 5.1. We  assert that the
 following implication holds: the   random variable $ \ \xi \ $ is such that

 $$
 \lim_{\lambda \to \infty} \phi^{-1}( \ln {\bf E} \exp(\lambda \ \xi))/\lambda = K = \const \in (0, \infty) \eqno(7.1)
 $$
if and only if

$$
\lim_{x \to \infty} \left( \phi^* \right)^{-1} ( |\ln {\bf P} (\xi \ge x)|)/x = 1/K. \eqno(7.2)
$$

 \vspace{4mm}

\ {\bf Remark 7.1.} This result is the direct generalization one represented in [11], theorem 1.4.1.  Wherein our proof is alike. \par

\vspace{4mm}

 \ {\bf 1.} We can and will suppose without loss of generality $ \ K = 1. \ $ \par
 Suppose at first  that the relation (7.2) be given.  Let $ \ \delta = \const \in (0, 1/e)  \ $ be arbitrary "small" number. Then there exists
a positive value $ \ x_0 = x_0(\delta) \  $   such that for all the values $ \ x \ge x_0(\delta) \ $ there holds

$$
 \left| (\phi^*)^{-1} (\ln T_{\xi}(x))/x - 1    \right| \le \delta.
$$
 \ We have solving the last inequality relative the tail function  $ \ T_{\xi}(x): \ $

$$
\exp \left[ \ - \phi^*(x(1+\delta)) \ \right] \le T_{\xi}(x) \le \exp \left[ \ - \phi^*(x(1-\delta)) \ \right], \eqno(7.3)
$$
and after integration by parts

$$
C_1 \lambda \int_R \exp \left[ \ \lambda x - \phi^*(x(1 + \delta)) \ \right] \ dx \le {\bf E}\exp(\lambda \ \xi) \le
$$

$$
C_2 \lambda \int_R \exp \left[ \ \lambda x - \phi^*(x(1 - \delta)) \ \right] dx, \ \lambda \ge 1. \eqno(7.4)
$$
 \ We conclude relaying  the proposition of theorem 2.1 for sufficiently greatest values of parameter  $ \ \lambda \ $

$$
\exp \phi(\lambda(1 - 2 \delta)) \le {\bf E}\exp(\lambda \ \xi) \le \exp \phi(\lambda(1 + 2 \delta)). \eqno(7.5)
$$

 \ By logarithm and taking the inverse function, we arrive at (7.1). \par

\vspace{4mm}

 \ {\bf 2.}  Conversely, suppose (7.1) holds true. Then for greatest positive  values $ \ \lambda \ $

$$
\exp \phi(\lambda(1 - \delta)) \le {\bf E} \exp(\lambda \xi) \le \exp \phi(\lambda(1 + \delta)).
$$
 \ Theorem 5.1 common with Tchebychev-Tchernov inequality  give us the following bilateral estimate

$$
\exp \left[ - \phi^*(x(1 + C_3 \delta)) \right] \le T_{\xi}(x)  \le \exp \left[ - \phi^*(x(1 - C_4 \delta)) \right]. \eqno(7.6)
$$
 \ The  relation (7.2) follows immediately from (7.6) after simple calculations. \par

 \vspace{4mm}

 \ {\it 2. Richter's approach.} \par

 \vspace{4mm}

 \ The proposition of theorem 5.1 is not true in the case when $ \ \delta = 0, \ $ i.e. when

 $$
 {\bf E} e^{\lambda \ \xi} = e^{\phi(\lambda)}, \  \lambda \ge 1,  \eqno(7.7)
 $$
 is the well-known Richter's case. We assume in the sequel that the function $ \  \phi(\cdot)  \ $ satisfies all the conditions of theorem 5.1. \par

 \ But one can apply the proof of theorem 4.1, where as before

 $$
 \phi_1(\lambda) = \phi_2(\lambda) = \phi(\lambda), \ T_{\xi}(x) \le \exp(-\phi^*(x)),
 $$

$$
\Delta = c_1 = \const > 0, \ \phi(\lambda) = \lambda \ x_0 - \phi^*(x_0).
$$

 \ We obtain after some calculations from (7.7)

$$
 T_{\xi}(x) \ge \exp(-\phi^*(x) - c_2[\phi] x), \ x \ge e.
$$

 \ To summarize: \par

 \vspace{4mm}

  \ {\bf Theorem 7.2.} We propose under formulated in this subsection conditions

$$
    \exp(-\phi^*(x) - c_2[\phi] x) \le  T_{\xi}(x) \le \exp(-\phi^*(x)), \ x \ge e. \eqno(7.8)
$$

\vspace{4mm}

\section{  \ Concluding remarks.}

 \vspace{3mm}

 \ It is interest by our opinion to generalize obtained in this article results onto a multidimensional random vector.  \par

 \vspace{6mm}

 {\bf References.}

 \vspace{4mm}


 \vspace{3mm}

{\bf 1.  Buldygin V.V., Kozachenko Yu.V. }  {\it Metric Characterization of Random
Variables and Random Processes.} 1998, Translations of Mathematics Monograph, AMS, v.188. \\

\vspace{3mm}

{\bf 2. Peter Eichelsbacher and Lucas Knichel.} {\it Fine asymptotics for models with Gamma type moments.} \\
arXiv:1710.06484v1 [math.PR] 17 Oct 2017 \\

 \vspace{3mm}

 {\bf 3. A. Fiorenza.}   {\it Duality and reflexivity in grand Lebesgue spaces. } Collect. Math.
{\bf 51,}  (2000), 131  - 148. \\

 \vspace{3mm}

{\bf  4. A. Fiorenza and G.E. Karadzhov.} {\it Grand and small Lebesgue spaces and
their analogs.} Consiglio Nationale Delle Ricerche, Instituto per le Applicazioni
del Calcoto Mauro Picone", Sezione di Napoli, Rapporto tecnico 272/03, (2005).\\

 \vspace{3mm}

{\bf 5. Alberto Fiorenza, Maria Rosaria Formica, Amiran Gogatishvili, Tengiz Kopaliani, Michel Rakotoson.}
{\it Characterization of interpolation between Grand, small or classical Lebesgue spaces. } \\
arXiv:1709.05892v1  [math.FA]  18 Sep 2017 \\

\vspace{3mm}

{\bf 6.  T. Iwaniec and C. Sbordone.} {\it On the integrability of the Jacobian under minimal
hypotheses. } Arch. Rat.Mech. Anal., 119, (1992), 129-143. \\

 \vspace{3mm}

{\bf 7. Janson.} {\it Further examples with moments of Gamma type.}\\
 arXiv:1204.56372v2, 2010.\\

\vspace{3mm}

{\bf 8. Janson.} {\it Moments of Gamma type and the Brownian supremum process area.}
 Probab. Surveys, {\bf 7,} (2010), 1–52.  MR 2645216

\vspace{3mm}

{\bf 9. Kozachenko Yu. V., Ostrovsky E.I. }  (1985). {\it The Banach Spaces of random Variables of subgaussian Type. } Theory of Probab.
and Math. Stat. (in Russian). Kiev, KSU, 32, 43-57. \\

\vspace{3mm}

{\bf 10. Kozachenko Yu.V., Ostrovsky E., Sirota L.}  {\it Relations between exponential tails, moments and
moment generating functions for random variables and vectors.} \\
arXiv:1701.01901v1 [math.FA] 8 Jan 2017 \\

 \vspace{3mm}

{\bf 11. Ostrovsky E.I. } (1999). {\it Exponential estimations for Random Fields and its
applications,} (in Russian). Moscow-Obninsk, OINPE. \\

 \vspace{3mm}

{\bf 12. Ostrovsky E. and Sirota L.} {\it Vector rearrangement invariant Banach spaces
of random variables with exponential decreasing tails of distributions.} \\
 arXiv:1510.04182v1 [math.PR] 14 Oct 2015 \\

 \vspace{3mm}

{\bf 13. Ostrovsky E. and Sirota L.}  {\it Non-asymptotical sharp exponential estimates
for maximum distribution of discontinuous random fields. } \\
 arXiv:1510.08945v1 [math.PR] 30 Oct 2015 \\

 \vspace{3mm}

{\bf 14. Ostrovsky E.I.}  {\it About supports of probability measures in separable Banach
spaces.} Soviet Math., Doklady, (1980), V. 255, $ \ N^0 \ $ 6, p. 836-838, (in Russian).\\

 \vspace{3mm}

\end{document}